\documentclass[12pt]{article}
\usepackage{amssymb}

\newcommand{\be}{\begin{equation}}
\newcommand{\ee}{\end{equation}}
\newcommand{\bea}{\begin{eqnarray}}
\newcommand{\eea}{\end{eqnarray}}
\newcommand{\barray}{\begin{array}}
\newcommand{\earray}{\end{array}}
\newcommand{\pa}{\partial}
\newcommand{\nn}{\nonumber}
\newcommand{\bitem}{\begin{itemize}}
\newcommand{\eitem}{\end{itemize}}
\newtheorem{teo}{Theorem}[section]
\newcommand{\bt}{\begin{teo}}
\newcommand{\et}{\end{teo}}
\newtheorem{Def}{Definition}[section]
\newcommand{\bd}{\begin{Def}}
\newcommand{\ed}{\end{Def}}
\newtheorem{lem}{Lemma}[section]
\newcommand{\bl}{\begin{lem}}
\newcommand{\el}{\end{lem}}
\newtheorem{prop}{Proposition}[section]
\newcommand{\bp}{\begin{prop}}
\newcommand{\ep}{\end{prop}}
\newtheorem{cor}{Corollary}[section]
\newcommand{\bc}{\begin{cor}}
\newcommand{\ec}{\end{cor}}
\newtheorem{ex}{Example}[section]
\newcommand{\bex}{\begin{ex}}
\newcommand{\eex}{\end{ex}}
\newtheorem{rem}{Remark}[section]
\newcommand{\br}{\begin{rem}}
\newcommand{\er}{\end{rem}}

\catcode `\@=11
\@addtoreset{equation}{section}

\begin{document}

\begin{center}
{\Large \textbf{The classification
of nonsingular multidimensional
Dubrovin--Novikov brackets\footnote{The research was supported by
the Max-Planck-Institut f\"{u}r Mathematik (Bonn, Germany), the
Russian Science Support Foundation, the Russian
Foundation for Basic Research
(grants Nos. 03-01-00782 and 05-01-00170), and the Program of Support for
Leading Scientific Schools (grant No. NSh-4182.2006.1).}}}
\end{center}

\medskip

\begin{center}
{\large {O. I. Mokhov}}
\end{center}
\begin{center}
{Centre for Nonlinear Studies,
L.D.Landau Institute for Theoretical Physics,\\
Russian Academy of Sciences,
Kosygina 2, Moscow, 117940, Russia\\
Department of Geometry and Topology,
Faculty of Mechanics and Mathematics, \\
M.V.Lomonosov Moscow State University,
Moscow, 119992, Russia\\
e-mail: mokhov@mi.ras.ru; mokhov@landau.ac.ru; mokhov@bk.ru}
\end{center}

\medskip

\begin{flushright}
{To the memory of my wonderful mother \hspace{13mm} \ \\
Maya Nikolayevna Mokhova
(04.05.1926 -- 12.09.2006)} \ \
\end{flushright}

\medskip

\begin{abstract}
In this paper the well-known Dubrovin--Novikov problem
posed as long ago as 1984 in connection with the Hamiltonian theory
of systems of hydrodynamic type, namely, the classification problem for
multidimensional Poisson brackets of hydrodynamic type, is solved.
In contrast to the one-dimensional case, in the general case,
a nondegenerate multidimensional Poisson bracket of hydrodynamic type
cannot be reduced to a constant form by a local change of coordinates.
Generally speaking, such a Poisson bracket is generated by
a nontrivial special infinite-dimensional Lie algebra.
In this paper we obtain a classification of all nonsingular
nondegenerate multidimensional Poisson brackets of hydrodynamic type for
any number $N$ of components and for any dimension $n$ by
differential-geometric methods. A key role
in the solution of this problem was
played by the theory of compatible metrics that had been earlier
constructed by the present author.
\end{abstract}



\section{Introduction}

In this paper the well-known Dubrovin--Novikov problem is solved.
This is the problem,
posed as long ago as 1984 in \cite{1}, of
the classification of
{\it multidimensional Poisson brackets
of hydrodynamic type}, i.e.,
field-theoretic Poisson brackets of the form
\be
\{ u^i (x), u^j (y) \} = \sum_{\alpha = 1}^n
\left ( g^{ij\alpha} (u(x)) \delta_{\alpha} (x - y)
+ b^{ij\alpha}_k (u(x)) u^k_{\alpha} (x) \delta (x - y)
\right ), \label{1}
\ee
where $u = (u^1, \ldots, u^N)$ are local coordinates on
a certain smooth $N$-dimensional manifold $M$ or in a
domain of ${\bf R}^N$;
$x = (x^1, \ldots, x^n),$ $y = (y^1, \ldots, y^n)$
are independent variables; the coefficients
$g^{ij\alpha} (u),$ $b^{ij\alpha}_k (u)$ are smooth
functions of the local coordinates $(u^1, \ldots, u^N)$,
$1 \leq i, j, k \leq N,$ $1 \leq \alpha \leq n$;
$u (x) = (u^1 (x), \ldots, u^N (x))$ are smooth functions
($N$-component fields) of
$n$ independent variables $x^1, \ldots, x^n$
with values in the manifold $M$; $u^k_{\alpha} (x) =
\pa u^k / \pa x^{\alpha}$; and
$\delta (x)$ is the Dirac delta-function,
$\delta_{\alpha} (x - y) =
\pa \delta (x - y) / \pa x^{\alpha}$.

The condition of skew-symmetry and the Jacobi identity
for a Poisson bracket of the form
(\ref{1}) impose very severe restrictions on the coefficients
$g^{ij\alpha} (u)$ and $b^{ij\alpha}_k (u)$ (we shall talk about these
restrictions below).

The general class of Poisson brackets of
the form (\ref{1}) was introduced by B.A.Dub\-rovin and S.P.Novikov
in \cite{2} (the one-dimensional case $n = 1$) and \cite{1}
(the multidimensional case of arbitrary $n$) in connection with
the Hamiltonian theory of systems of hydrodynamic type,
and such Poisson brackets are also called by
{\it multidimensional Dubrovin--Novikov brackets}.
The very natural differential-geometric Hamiltonian approach
to {\it systems of hydrodynamic type}, i.e.,
evolution systems of homogeneous quasilinear equations with
partial derivatives of the first order, proposed in
\cite{1} and \cite{2}, is connected with the general
class of Poisson brackets of the form (\ref{1}).

If $\det (g^{ij\alpha} (u)) \neq 0$ for all $\alpha$,
then the bracket (\ref{1}) is called {\it nondegenerate}.
The one-dimensional nondegenerate Poisson brackets of hydrodynamic type
were completely classified
by B.A.Dubrovin and S.P.Novikov in 1983 in
\cite{2}, where they proved that each such Poisson bracket
can be reduced to a constant form by a local change of coordinates
and can be characterized by a unique
invariant, namely, by signature of the
metric of the bracket.

Afterwards this Hamiltonian approach led to the creation of the remarkable
and very fruitful theory of integrable one-dimensional
systems of hydrodynamic type (see \cite{2a}, \cite{2b}),
and at present one can consider that the Hamiltonian
and integrability properties of one-dimensional
systems of hydrodynamic type are well studied
(although in the one-dimensional case also many
very important problems remain still unsolved).

However, the case of multidimensional
systems of hydrodynamic type is considerably more complicated and,
in fact, has not been studied as yet from the viewpoint of
the general Hamiltonian
and integrability properties
(see \cite{2c}).
In the multidimensional case, very serious problems
arise already under studying the simplest natural invariant
class of suitable local Hamiltonian structures
(\ref{1}).

In \cite{1} B.A.Dubrovin and S.P.Novikov
posed the classification problem for nondegenerate
multidimensional Poisson brackets of the form
(\ref{1}) and they showed therewith that for the
multidimensional brackets,
in contrast to one-dimensional, the statement on reducing
to a constant form by a local change of coordinates is
definitely not true, since, in particular,
the two-dimensional nondegenerate Poisson bracket of
hydrodynamic type, generated by the Lie algebra of vector fields
on the plane and related to the two-dimensional Euler hydrodynamics,
cannot be reduced to a constant form by a local change of coordinates
(there are nonzero tensor obstructions in this case).
In the same paper \cite{1} B.A.Dubrovin and S.P.Novikov
undertook an unsuccessful attempt to obtain a classification
of these brackets
in the case $N=1$ and in the first nontrivial case
$n = N = 2$ (it was shown in \cite{3} by the author of the present paper
that these results in the paper \cite{1} were erroneous).

A complete classification of all one-component
($N=1$, $n$ arbitrary) and all nondegenerate two-component
($N = 2$, $n$ arbitrary) Poisson brackets of the form
(\ref{1}) was obtained by the present author in \cite{3},
where a complete classification
was also obtained for
$1 \leq N \leq 4$, for which it was necessary
to develop powerful algebraic machinery, but all attempts to solve
the problem for an arbitrary number
$N$ of components by algebraic methods have so far been
unsuccessful.

In the general case, a nondegenerate multidimensional
Poisson bracket of hydrodynamic type cannot
be reduced to a constant form by a local change of
coordinates, in contrast to the one-dimensional case,
but according to results of B.A.Dubrovin and S.P.Novikov \cite{1},
and also the present author
\cite{3}, for any such Poisson bracket,
there always exist special local coordinates
such that the bracket is linear with respect to the fields
(generally speaking, linear nonhomogeneous).

Multidimensional
Poisson brackets of hydrodynamic type
that are linear with respect to the fields are generated
by special infinite-dimensional Lie algebras
({\it Lie algebras of hydrodynamic type})
and special 2-cocycles on these
infinite-dimensional Lie algebras (for nonhomogeneous
linear brackets). A classification of these Lie algebras for
$n \geq 2$ was obtained
for the cases
$1 \leq N \leq 4$ by the author of the present paper in \cite{3}.
(One-dimensional Dubrovin--Novikov brackets
that are linear with respect to the fields and
algebraic structures related to them were studied
in \cite{4}.)

The classification problem for
multidimensional nondegenerate
Poisson brackets of hydrodynamic type, for an arbitrary
number of components, was absolutely inaccessible and inapproachable
up to now. In this paper we obtain a classification of all
nonsingular nondegenerate multidimensional
Poisson brackets of the form
(\ref{1}) (in any case, for the situation
of brackets of the form (\ref{1}) in general position)
for arbitrary number $N$ of components and for arbitrary
$n$ by differential-geometric methods.
A key role
in the solution of this problem was
played by the theory of compatible metrics
constructed by the present author in \cite{5}.

\section{Multidimensional homogeneous
brackets of\\ the first order and Dubrovin--Novikov brackets}

The requirement of bilinearity and the fulfilment of
the Leibniz identity for a bracket of the form
(\ref{1}), i.e., an arbitrary {\it multidimensional
homogeneous bracket of the first order}, is equivalent to the condition
that on arbitrary functionals $I$ and $J$ on the space of fields
$u(x)$ the bracket has the form
\be
\{ I, J \} = \sum_{\alpha = 1}^n \int {\delta I \over \delta u^i (x)}
\left ( g^{ij\alpha} (u(x)) {d \over d x^{\alpha}}
+ b^{ij\alpha}_k (u(x)) u^k_{\alpha} (x)
\right ) {\delta J \over \delta u^j (x)} d^n x, \label{1n}
\ee
where $d / d x^{\alpha}$ is the total derivative
with respect to the independent variable $x^{\alpha}$,
\bea
&&
{d \over d x^{\alpha}} =
{\pa \over \pa x^{\alpha}} + u^i_{\alpha}
{\pa \over \pa u^i} + u^i_{\beta \alpha}
{\pa \over \pa u^i_{\beta}} + \cdots +
u^i_{\beta_1 \cdots \beta_s \alpha}
{\pa \over \pa u^i_{\beta_1 \cdots \beta_s}} + \cdots, \\
&& \hspace{30mm}
u^i_{\beta_1 \cdots \beta_s} =
{\pa^s u^i \over \pa x^{\beta_1} \cdots \pa x^{\beta_s}}, \nn
\eea
where summation over repeating upper and lower indices is assumed,
moreover, in this case the summation over $\beta_1, \ldots, \beta_s$
is taken over distinct (up to arbitrary permutations)
sets of these indices (one can consider that these sets of indices are
ordered:
$\beta_1  \leq \cdots \leq \beta_s$).

The class of brackets of the form (\ref{1n}) is
invariant with respect to local changes of coordinates
$(u^1, \ldots, u^N)$.  Here for each $\alpha$ the coefficients
$g^{ij\alpha} (u)$ and $b^{ij\alpha}_k (u)$ of the brackets
are transformed as differential-geometric objects.
For any bracket of the form
(\ref{1n}) the coefficient
$g^{ij\alpha} (u)$ for each $\alpha$ is a contravariant two-valent
tensor on the manifold $M$. For any nondegenerate bracket
of the form (\ref{1n}) and for each $\alpha$ we can introduce
coefficients
$\Gamma^{k\alpha}_{ij} (u)$ by the formulae
$\Gamma^{k\alpha}_{ij} (u) = -g_{is}^\alpha (u) b^{sk\alpha}_j (u)$,
where $g_{is}^\alpha (u)$ is a covariant two-valent
tensor inverse to the tensor $g^{ij\alpha} (u)$,
$g^{ij\alpha} (u)g_{jk}^\alpha (u) = \delta^i_k$.
Here for any nondegenerate bracket of the form
(\ref{1n}) and for each $\alpha$ the coefficients
$\Gamma^{j\alpha}_{sk} (u)$ are transformed
as the coefficients of an affine connection
on the manifold $M$. Thus a nondegenerate
bracket of the form (\ref{1n}) is given by $n$
arbitrary nondegenerate contravariant
two-valent tensors $g^{ij\alpha} (u)$, $1 \leq \alpha \leq n$,
and $n$ arbitrary affine connections
$\Gamma^{j\alpha}_{sk} (u)$,
$1 \leq \alpha \leq n$, on the manifold $M$.

Note that among brackets of the form (\ref{1n})
there are, in particular, very simple, but very important
brackets that do not depend on the fields
$u(x)$ at all, namely, {\it constant brackets}
({\it brackets of a constant form}). All such brackets
are given by the conditions $g^{ij\alpha} (u) = {\rm \ const}$,
$b^{ij\alpha}_k (u) = 0$.
The problem, which we solve in this paper,
is to find what the simplest canonical forms can be reduced
Poisson brackets of the form (\ref{1n})
by local changes of coordinates to.

First of all, we shall show the general relations
on the coefficients of Poisson brackets of the form
(\ref{1n}) \cite{6}.

\bt  {\rm \cite{6}}. \label{t1}
A bracket {\rm (\ref{1n})} is a Poisson bracket, i.e.,
is skew-symmetric and satisfies the Jacobi identity,
if and only if the following relations for the coefficients of
the bracket are fulfilled{\rm :}
\be
g^{ij\alpha}=g^{ji\alpha},
\label{a1}
\ee
\be
\frac{\pa g^{ij\alpha}}{\pa u^k}=b^{ij\alpha}_k+b^{ji\alpha}_k,
\label{a2}
\ee
\be
\sum_{(\alpha,\beta)}\left (g^{si\alpha}b^{jr\beta}_s
-g^{sj\beta}b^{ir\alpha}_s\right )=0,
\label{a3}
\ee
\be
\sum_{(i,j,r)}\left (g^{si\alpha}b^{jr\beta}_s
-g^{sj\beta}b^{ir\alpha}_s\right )=0,
\label{a4}
\ee
\be
\sum_{(\alpha,\beta)}\left [g^{si\alpha}
\left (\frac{\pa b^{jr\beta}_s}{\pa u^q}
-\frac{\pa b^{jr\beta}_q}{\pa u^s}\right )
+b^{ij\alpha}_sb^{sr\beta}_q-b^{ir\alpha}_sb^{sj\beta}_q\right ]=0,
\label{a5}
\ee
\be
g^{si\beta}\frac{\pa b^{jr\alpha}_q}{\pa u^s}
-b^{ij\beta}_sb^{sr\alpha}_q-b^{ir\beta}_sb^{js\alpha}_q =
g^{sj\alpha}\frac{\pa b^{ir\beta}_q}{\pa u^s}
-b^{ji\alpha}_sb^{sr\beta}_q-b^{is\beta}_qb^{jr\alpha}_s,
\label{a6}
\ee
\bea
&&
\frac\pa{\pa u^k}
\left [g^{si\alpha}\left (\frac{\pa b^{jr\beta}_s}{\pa u^q}
-\frac{\pa b^{jr\beta}_q}{\pa u^s}\right )
+b^{ij\alpha}_sb^{sr\beta}_q-b^{ir\alpha}_sb^{sj\beta}_q\right ] +
\nn\\
&&
+\sum_{(i,j,r)}\left [b^{si\beta}_q
\left (\frac{\pa b^{jr\alpha}_k}{\pa u^s}
-\frac{\pa b^{jr\alpha}_s}{\pa u^k}\right )\right ] + \nn\\
&&
+\frac\pa{\pa u^q}
\left [g^{si\beta}\left (\frac{\pa b^{jr\alpha}_s}{\pa u^k}
-\frac{\pa b^{jr\alpha}_k}{\pa u^s}\right )
+b^{ij\beta}_sb^{sr\alpha}_k-b^{ir\beta}_sb^{sj\alpha}_k\right ]+
\nn\\
&&
+\sum_{(i,j,r)}\left [b^{si\alpha}_k
\left (\frac{\pa b^{jr\beta}_q}{\pa u^s}
-\frac{\pa b^{jr\beta}_s}{\pa u^q}\right )\right ]=0.
\label{a7}
\eea
Relations {\rm (\ref{a1})} and {\rm (\ref{a2})}
are equivalent to the skew-symmetry of a bracket
{\rm (\ref{1n})}, and relations {\rm (\ref{a3}){--}(\ref{a7})}
are equivalent to the fulfilment of the Jacobi identity
for a skew-symmetric bracket of the form {\rm (\ref{1n})}.
\et
We do not assume nondegeneracy of brackets and do not
impose any additional conditions
on the coefficients of brackets (\ref{1n})
in Theorem \ref{t1}.
Moreover, this theorem remains true for the important case
(for example, in the theory of nonlinear chains), when  the indices
run an infinite set of values (for example,
the sets of all integer or natural numbers) and any function
under consideration (i.e., all the coefficients of brackets)
depends on an arbitrary, but finite, number of the variables
(in the one-dimensional case, a similar observation  was made
by I.Dorfman in \cite{6a}).
The signs $\sum_{(\alpha, \beta)}$ and
$\sum_{(i, j, k)}$ mean summations over all
cyclic permutations of the indicated indices, i.e.,
in the given case, the indices
$(\alpha, \beta)$ and $(i, j, k)$ respectively.
In the one-dimensional case the general relations
on the coefficients of the Poisson brackets of hydrodynamic type
(without the assumption of nondegeneracy)
were shown in \cite{7}, where there is an error in the formulae,
corrected in \cite{8} under construction of the nonlocal
generalization of the one-dimensional
Dubrovin--Novikov brackets.

In particular, the following important lemma immediately follows from
the relations of Theorem \ref{t1}.

\bl \label{l1}
For each multidimensional Poisson bracket of the form
{\rm (\ref{1n})} and for each $\alpha$ the corresponding
summand on the right-hand side of the formula {\rm (\ref{1n})}
is a one-dimensional Poisson bracket of hydrodynamic type,
i.e., each multidimensional Poisson bracket of the form
{\rm (\ref{1n})} is always the sum of one-dimensional
Poisson brackets with respect to each of the independent variables
$x^\alpha$.
\el

Besides, in the multidimensional case, for
$n \geq 2$, the relations of Theorem \ref{t1}
impose additional nontrivial restrictions on the coefficients
of these one-dimensional Poisson brackets.

First of all, let us consider the relations
(\ref{a1})--(\ref{a7}) in the one-dimensional case
(for $n=1$). It is obvious that in the one-dimensional
case the relation (\ref{a4}) follows from the relation
(\ref{a3}) and, moreover, in this case
the relation (\ref{a6}) follows from the relations
(\ref{a5}), (\ref{a3}) and (\ref{a2}).
In the nondegenerate one-dimensional case
($\det (g^{ij\alpha} (u)) \neq 0$, $x = x^{\alpha}$) we have:
the relation (\ref{a1}) gives the condition that
a nondegenerate tensor
$g^{ij\alpha}$ is symmetric, i.e.,
$g^{ij\alpha}$ is a Riemannian or pseudo-Riemannian
contravariant metric on the manifold $M$;
the relation (\ref{a2}) means that the connection
$\Gamma^{k\alpha}_{ij}=-g_{is}^\alpha b^{sk\alpha}_j$
is compatible with the metric $g^{ij\alpha}$, i.e.,
the corresponding covariant derivative of the metric is equal to zero;
the relation (\ref{a3}) is equivalent to
the condition that the connection $\Gamma^{k\alpha}_{ij}$
is symmetric, i.e.,
$\Gamma^{k\alpha}_{ij}=\Gamma^{k\alpha}_{ji}$;
the relation (\ref{a5}) means exactly that
the connection is flat, i.e.,
the Riemannian curvature tensor vanishes;
the relation (\ref{a7}) for
$\alpha=\beta$ and for a nondegenerate metric $g^{ij\alpha}$
follows from the relations (\ref{a3}) and (\ref{a5}).

This proves the Dubrovin--Novikov theorem \cite{2}:
an arbitrary nondegenerate one-dimensional Poisson bracket of
hydrodynamic type is uniquely determined by an arbitrary
flat metric $g^{ij} (u)$ and is reduced to
a constant form in any flat coordinates of this metric,
and $\Gamma^k_{ij} (u) =-g_{is} (u) b^{sk}_j (u)$ is the
Levi-Civita connection generated by the metric $g^{ij} (u)$.

Thus, by virtue of Lemma \ref{l1} and the Dubrovin--Novikov theorem
for nondegenerate one-dimensional Poisson brackets of the
form (\ref{1}) \cite{2}, all the tensors
$g^{ij\alpha} (u)$ are flat metrics (metrics of
zero Riemannian curvature), and
each the corresponding affine connection $\Gamma^{j\alpha}_{sk} (u)$
is compatible with the respective metric $g^{ij\alpha} (u)$
and has zero torsion and zero Riemannian curvature,
i.e., the affine connection $\Gamma^{j\alpha}_{sk} (u)$ is
the Levi-Civita connection and is uniquely determined
by the flat metric $g^{ij\alpha} (u)$.

Hence,
each nondegenerate Poisson bracket of the form
(\ref{1n}) is uniquely determined by the flat metrics
$g^{ij\alpha} (u)$, which are connected by
additional severe restrictions in the multidimensional case,
and our problem is reduced to a classification of the admissible
sets of the
flat metrics $g^{ij\alpha} (u)$, $1 \leq \alpha \leq n$.

\section{Tensor obstructions and tensor relations}

Let us consider an arbitrary nondegenerate
multidimensional Poisson bracket of the form
(\ref{1n}) and introduce the tensors
$T^{i\alpha\beta}_{jk} (u) =
\Gamma^{i\beta}_{jk} (u) - \Gamma^{i\alpha}_{jk} (u)$
defined for each pair of distinct indices
$\alpha$ and $\beta$.
The tensors $T^{i\alpha\beta}_{jk} (u)$
are obstructions for reducing nondegenerate
multidimensional Poisson brackets of the form
(\ref{1n}) to constant brackets, i.e.,
an arbitrary nondegenerate Poisson bracket
of the form (\ref{1n}) can be reduced
to a constant bracket by a local change
of coordinates if and only if
all the obstruction tensors $T^{i\alpha\beta}_{jk} (u)$
are identically equal to zero.

Indeed, if even one of these tensors is not
equal to zero identically, then
the Poisson bracket cannot be reduced to
a constant form by a local change
of coordinates, since for any constant bracket
the coefficients of all the connections $\Gamma^{i\alpha}_{jk}(u)$
are equal to zero identically, and, consequently,
in these coordinates all the tensors
$T^{i\alpha\beta}_{jk}(u)$ should be equal to zero identically.
The converse statement is also obvious. Indeed,
if all the tensors
$T^{i\alpha\beta}_{jk}(u)$ are equal to zero identically, then
all the connections $\Gamma^{i\alpha}_{jk}(u)$
are equal one to another and, consequently,
all of them are equal to zero in any flat coordinates
of the metric $g^{ij1}(u)$, and all the metrics of the bracket
are necessarily constant in these coordinates
by virtue of compatibility of the metrics with the
corresponding connections.

The following theorem gives a complete set of
tensor relations defining the class of
nondegenerate Poisson brackets of the form (\ref{1n}).

\bt {\rm \cite{3}}. \label{t2}
Flat nondegenerate metrics $g^{ij\alpha} (u)$
define a multidimensional Poisson bracket of the form
{\rm (\ref{1n})} if and only if the following
relations are fulfilled{\rm :}
\be
T^{ijk\alpha\beta} (u) =
T^{kji\alpha\beta} (u), \label{b1}
\ee
\be
\sum_{(i, j, k)} T^{ijk\alpha\beta} (u) = 0, \label{b2}
\ee
\be
T^{ijs\alpha\beta} (u) T^{r\alpha\beta}_{st} (u)
= T^{irs\alpha\beta} (u) T^{j\alpha\beta}_{st} (u), \label{b3}
\ee
\be
\nabla^{\alpha}_r T^{ijk\alpha\beta} (u) = 0, \label{b4}
\ee
where
$T^{i\alpha\beta}_{jk} (u) =
\Gamma^{i\beta}_{jk} (u) - \Gamma^{i\alpha}_{jk} (u)$,
$T^{ijk\alpha\beta} (u) =
g^{ks\beta} (u) g^{ir\alpha} (u) T^{j\alpha\beta}_{rs} (u)$,
the sign $\sum_{(i, j, k)}$ means summation
over all cyclic permutations of indices $(i, j, k)$,
$\nabla^{\alpha}_r$ is the covariant derivative
given by the connection $\Gamma^{i\alpha}_{jk} (u)$,
and $\Gamma^{i\alpha}_{jk} (u)$ is the Levi-Civita connection
generated by the metric $g^{ij\alpha} (u)$.
\et

The relations (\ref{b1}) and (\ref{b3}) were found
by B.A.Dubrovin and S.P.Novikov in \cite{1}.
It was proved by the present author in
\cite{3} that this set of tensor relations is not complete,
and the complete system of tensor relations
(\ref{b1})--(\ref{b4}) for nondegenerate Poisson brackets of the
form (\ref{1n}) was obtained \cite{3}.
The relations (\ref{b2}) and
(\ref{b4}) are very essential and play an important role
under classifying the multidimensional Poisson brackets of
hydrodynamic type.

If all the metrics
$g^{ij\alpha}$ are nondegenerate, then we obtain
from Theorem \ref{t1} that
for $\alpha\ne\beta$ the condition (\ref{a3})
gives the relation (\ref{b1}) for the obstruction tensors;
the condition (\ref{a4}) is equivalent to the relation
(\ref{b2}); the condition (\ref{a5}) is equivalent
to the relation (\ref{b3}); the condition
(\ref{a6}) is equivalent to the relation (\ref{b4});
the condition (\ref{a7}) for nondegenerate metrics
is a direct consequence of the relations
(\ref{a1})--(\ref{a6})
(this is not true in the general case and so the condition
(\ref{a7}) is essential in the case of degenerate metrics).

It is not complicated to obtain a complete
classification of all one-component (degenerate and nondegenerate)
Poisson brackets
of the form (\ref{1n}).
Let us prove that for $N=1$ and arbitrary $n$
all the obstruction tensors are identically equal to zero.
Indeed, in the one-component case the relations
(\ref{a1}),
(\ref{a3}), (\ref{a5}) and (\ref{a7}) are automatically fulfilled,
the relation (\ref{a2}) gives
$\pa g^\alpha/\pa u=2b^\alpha (u)$, from the relation
(\ref{a4}) we obtain that
$g^\alpha (u) b^\beta (u) = g^\beta (u) b^\alpha (u)$, and
the relation
(\ref{a6}) follows from (\ref{a2}) and (\ref{a4}).
Thus in the one-component case for any $n$ and
for any indices $\alpha$ and
$\beta$ we have
$\Gamma^{\alpha} (u) = - b^{\alpha} (u)/ g^{\alpha} (u) =
 - b^{\beta} (u)/ g^{\beta} (u)
= - (1/(2g^{\alpha}(u))) \pa g^{\alpha} / \pa u
= - (1/(2g^{\beta}(u))) \pa g^{\beta} / \pa u
= \Gamma^{\beta} (u)$, i.e.,
 $g^\alpha (u)=c^\alpha g(u)$, where $g (u)$ is
an arbitrary nonzero function, and $c^\alpha$ are arbitrary
nonzero constants. Hence, all the obstruction
tensors $T^{\alpha\beta} (u) = \Gamma^{\beta} (u) -
\Gamma^{\alpha} (u)$ are identically equal to zero,
and each multidimensional one-component
Poisson bracket of the form (\ref{1n})
is reduced to a constant form by a local change
of the unique coordinate $u = u^1$ (in the degenerate case,
i.e., if $g^\alpha (u)=0$, then also $b^\alpha (u)=0$
by virtue of the relation (\ref{a2})). Hence, we obtain
a complete classification of
all one-component (degenerate and nondegenerate)
Poisson brackets
of the form (\ref{1n}).

\section{Linear Poisson brackets of hydrodynamic type, \\
infinite-dimensional Lie algebras, Frobenius and \\
quasi-Frobenius algebras}

Each flat metric $g^{ij\alpha}(u)$ can be reduced
to a constant form by a local change of coordinates.
Let us reduce one of the metrics
(for definiteness, we reduce the first one)
of an arbitrary nondegenerate Poisson bracket of the form
(\ref{1n}) to a constant form.

\bt  \label{t3}
If for a nondegenerate multidimensional
Poisson bracket of the form {\rm (\ref{1n})}
we have $g^{ij1}={\rm const}$, then
either this Poisson bracket is constant in these local
coordinates, or it has nonzero obstruction tensors, i.e.,
this Poisson bracket cannot be reduced to
a constant form by local changes
of coordinates at all.
\et

Indeed, if in the given coordinates
$\Gamma ^{i\alpha}_{jk} = 0$ for all $\alpha$,
then $b^{ij\alpha}_k = 0$ for all $\alpha$ and the
Poisson bracket is constant by virtue of
the relation (\ref{a2}). But if for a certain $\alpha$
not all of the coefficients of the connection
$\Gamma ^{i\alpha}_{jk}$ are equal to zero,
then also the obstruction tensor
$T^{i1\alpha}_{jk}$, which coincides in the given local coordinates
with $\Gamma ^{i\alpha}_{jk}$, is not equal to zero,
and, consequently, this Poisson bracket cannot
be reduced to a constant form by local changes of coordinates.

The following important theorem is a simple
consequence of Theorem \ref{t2}.

\bt {\rm \cite{1}, \cite{3}}. \label{t4}
If $g^{ij1}={\rm const}$ for a nondegenerate
multidimensional Poisson bracket of the form
{\rm (\ref{1n})}, then all the other metrics
are linear with respect to the local coordinates
$u^i$, and the Poisson bracket is linear with respect to the fields
$u (x)${\rm :}
$$
g^{ij\alpha} (u) =
(b^{ij\alpha}_k+b^{ji\alpha}_k)u^k+g^{ij\alpha}_0,\qquad
b^{ij\alpha}_k={\rm const},\quad g^{ij\alpha}_0={\rm const},
\quad 2 \leq \alpha \leq n.
$$
\et

Let the metric $g^{ij1}$ be constant. Then in these
local coordinates we have
$\Gamma^{i1}_{jk} =0.$
It follows from the relation (\ref{b4}) that
in these coordinates
$T^{ijk1\alpha} = {\rm const}$ for all $\alpha$.
Let us prove that all the coefficients $b^{ij\alpha}_k$ are also
constant in these coordinates.
Indeed, we have
$$
b^{ij\alpha}_k = -g^{is\alpha} \Gamma^{j\alpha}_{sk}=
-g^{is\alpha} \Gamma^{j \alpha}_{ks}=
-g^{is\alpha} T^{j1 \alpha}_{ks}
= - g_{kq1}g^{is\alpha} g^{qr1} T^{j1\alpha}_{rs}=
- g_{kq1}T^{qji1\alpha} = {\rm const}.
$$
Here, we have used the fact that all the connections
$\Gamma ^{i\alpha}_{jk}$ are symmetric.
Then the linearity of all the metrics in these local coordinates
and the linearity of the bracket with respect to
the fields follows from the relation (\ref{a2}).

For $N\ge3$ Theorem \ref{t4} was proved by B.A.Dubrovin
and S.P.Novikov in \cite{1}. The study of the cases
$N = 1$ and $N = 2$ in \cite{1} is erroneous, since
the study is based on an incomplete set of relations for
the obstruction tensors, which is obtained in
\cite{1} and which is insufficient
in order to guarantee that the bracket
(\ref{1n}) is a Poisson bracket.
In the complete form Theorem \ref{t4} was proved by
the present author in \cite{3}.

Let us recall here very briefly, in a necessary for us form,
a general scheme (see also \cite{9}) that goes back to Sophus Lie
and concerns interconnections between Lie algebras and
Poisson structures whose coefficients depend linearly
(possibly, nonhomogeneously) on coordinates
(the Lie--Poisson brackets).

For the general infinite-dimensional case
we shall describe special infinite-dimen\-sional
Lie algebras corresponding to arbitrary Poisson structures
whose coefficients depend linearly
(possibly, nonhomogeneously) on the field variables
$u^i(x)$ and their derivatives
$u^i_{(k)}$, $(k)=(k_1,\dots,k_n)$, where
$u^i_{(k)}=
\pa^{|k|}u^i/(\pa(x^1)^{k_1}\cdots \pa(x^n)^{k_n})$,
$|k|=k_1+\cdots+k_n$.

An operator $M^{ij}$ determining a Poisson bracket
(a {\it Poisson structure})
\be
\{ I, J \} = \int {\delta I \over \delta u^i (x)}
M^{ij} {\delta J \over \delta u^j (x)} d^n x
\ee
is called {\it Hamiltonian}.

Consider arbitrary Hamiltonian operators
given by differential operators whose coefficients
depend linearly
(generally speaking, nonhomogeneously) on the
field variables $u^i(x)$ and their derivatives, i.e.,
Hamiltonian operators of the form
\be
M^{ij}=\left (a^{ij,(k)(p)}_su^s_{(k)}+b^{ij,(p)}\right )
\frac{d^{|p|}}{d(x^1)^{p_1}\cdots d(x^n)^{p_n}}\,,
\label{c1}
\ee
where $a^{ij,(k)(p)}_s$ and $b^{ij,(p)}$ are constants.

Consider the infinite-dimensional space $S$ of
sequences $(\xi_1,\dots,\xi_N)$, where $\xi_i\in C^\infty(T^n)$ are
smooth functions on an $n$-dimensional torus.
Then if
$\xi$ and $\eta$ belong to the space $S$, we have
\be
(\xi,M(\eta))\equiv\int_{T^n}\xi_iM^{ij}\eta_j\,d^nx
=\int_{T^n}u^s[\xi,\eta]_s\,d^nx+
\int_{T^n}\xi_ib^{ij,(p)}\eta_{j(p)}\,d^nx. \label{c0}
\ee
Thus on the space $S$ a bilinear operation
 $[\,\cdot\,{,}\,\cdot\,]$,
\be
(\xi,\eta)\mapsto\zeta=[\xi,\eta]\in S,
\ \ \ \ \zeta_s=[\xi,\eta]_s= \sum_{i, j, (k), (p)}
(- 1)^{|k|} a^{ij,(k)(p)}_s (\xi_i \eta_{j (p)})_{(k)}, \label{c2}
\ee
and a bilinear form
\be
\omega(\xi,\eta)=\int_{T^n}\xi_ib^{ij,(p)}\eta_{j(p)}\,d^nx \label{c3}
\ee
are defined.

An arbitrary operator $M^{ij}$ of the form (\ref{c1})
is skew-symmetric if and only if the corresponding
bilinear operation
(\ref{c2}) and bilinear form (\ref{c3}) are skew-symmetric on the space
$S$, i.e.,
$[\xi,\eta]=-[\eta,\xi]$
and
$\omega(\xi,\eta)=-\omega(\eta,\xi).$

An arbitrary operator $M^{ij}$ of the form (\ref{c1})
is Hamiltonian if and only if the corresponding space $S$
is a Lie algebra with respect to
the corresponding bilinear operation (\ref{c2}), i.e.,
this bilinear operation is skew-symmetric and satisfies
the Jacobi identity
$[\xi,[\eta,\zeta]]+[\eta,[\zeta,\xi]]+[\zeta,[\xi,\eta]]=0,$
and, in addition, the corresponding
bilinear form (\ref{c3}) is a $2$-cocycle
on this Lie algebra, i.e., this bilinear form is skew-symmetric
and satisfies the closedness identity
\be
(d\omega)(\xi,\eta,\zeta)\equiv
\omega([\xi,\eta],\zeta)+\omega([\eta,\zeta],\xi)+
\omega([\zeta,\xi],\eta)=0.\nn
\ee

Note that a $2$-cocycle $\omega(\xi,\eta)$ defined by
a Hamiltonian operator of the form
(\ref{c1}) is
{\it cohomologous to zero}, i.e.,
$\omega(\xi,\eta)=(df)(\xi,\eta)\equiv f([\xi,\eta]),$
where $f$ is a $1$-form on the Lie algebra $S$,
if and only if this $2$-cocycle can be annihilated by a shift of the
field variables
$u^i\mapsto u^i-c^i$, where $c^i$ are arbitrary constants,
i.e., provided that
$b^{ij,(p)}=a^{ij,(0)(p)}_kc^k.$

\bl \label{l2}
A multidimensional Poisson bracket of the form
{\rm (\ref{1n})} is linear
{\rm (}possibly, nonhomogeneously{\rm )}
with respect to the fields $u (x)$ if and only if
\be
g^{ij\alpha} (u) =
(b^{ij\alpha}_k+b^{ji\alpha}_k)u^k+g^{ij\alpha}_0,\quad
b^{ij\alpha}_k={\rm const},\quad g^{ij\alpha}_0={\rm const},
\quad 1 \leq \alpha \leq n, \label{m1}
\ee
here the constants
$b^{ij\alpha}_k$ and $g^{ij\alpha}_0$ satisfy
nontrivial quadratic relations following from
the relations
{\rm (\ref{a1})}--{\rm (\ref{a7})} of Theorem
{\rm \ref{t1}} and define $n$
compatible algebras $ {\cal B^{\alpha}}$, $1 \leq \alpha \leq n,$
of Frobenius and quasi-Frobenius type with
ctructural constants $b^{ij\alpha}_k$ and
symmetric bilinear forms $g^{ij\alpha}_0$ on these algebras
{\rm (}one can consider that there is one $N$-dimensional
algebra with the basis
$e^1, \ldots, e^N$ furnished with
$n$ compatible multiplications and $n$
symmetric bilinear forms,
$e^i \stackrel{\alpha}{\circ} e^j = b^{ij\alpha}_k e^k$,
$\langle e^i, e^j \rangle_{\alpha} = g^{ij\alpha}_0$,
$1 \leq \alpha \leq n${\rm )}.
To each multidimensional Poisson bracket
of the form {\rm (\ref{1n})} that is linear
{\rm (}possibly, nonhomogeneously{\rm )} with respect to
the fields $u (x)$ it corresponds an infinite-dimensional Lie algebra
of a special type with a $2$-cocycle of a special type
on it {\rm (a Lie algebra of hydrodynamic type):}
\be
[\xi,\eta]_k=b^{ij\alpha}_k((\eta_i)_\alpha
\xi_j-\eta_j(\xi_i)_\alpha), \label{alx}
\ee
\be
b^{ij\alpha}_k = {\rm \ const},\quad
\xi=(\xi_1,\dots,\xi_N),\quad \xi_i(x)\in C^1(T^n), \label{al}
\ee
\be
\omega(\xi,\eta)=\int g^{ij\alpha}_0(\eta_j(x))_\alpha\xi_i(x)\,d^nx,
\quad g^{ij\alpha}_0 = {\rm \ const}. \label{ko}
\ee
\el

The classification of multidimensional Poisson brackets
of the form (\ref{1n}) that are linear
(possibly, nonhomogeneously) with respect to the fields
$u (x)$ and the classification of related to them
infinite-dimensional Lie algebras,
admissible $2$-cocycles on these Lie algebras, and
also algebras
${\cal B^{\alpha}}$, $1 \leq \alpha \leq n,$
of Frobenius and quasi-Frobenius type generating the brackets
is a separate important problem, which is far from a complete
solution for now.

In the one-dimensional case, this problem
was solved by A.A.Balinsky and S.P.No\-vikov in
\cite{4}, where, as is important to note,
was discovered its connection with the theory of Frobenius
and quasi-Frobenius algebras, although in the
one-dimensional case also many questions remained unsolved.

In the multidimensional case for
$n \geq 2$ this problem was studied by the present author in
\cite{3}, where the theory of the corresponding
algebras was developed and complete classification results
were obtained for the cases, when the number of
components is not more than 4.
The results on this problem and, in particular,
on the theory of the corresponding quasi-Frobenius algebras
will be published in a separate paper.

\bt \label{t5}
Each nondegenerate multidimensional
Poisson bracket of the form
{\rm (\ref{1n})} is defined by a certain
infinite-dimensional Lie algebra of the form
{\rm (\ref{alx}), (\ref{al})} with a certain
$2$-cocycle of the form {\rm (\ref{ko})} on this Lie algebra
for which
\be
b^{ij1}_k = 0,
\ \ \ \ \det \left ((b^{ij\alpha}_k+b^{ji\alpha}_k)u^k+
g^{ij\alpha}_0\right ) \neq 0, \quad 1 \leq \alpha \leq n. \label{m2}
\ee
\et
Theorem \ref{t5} immediately follows from
Theorem \ref{t4} and Lemma \ref{l2}.
The theory of infinite-dimensional Lie algebras
of hydrodynamic type (\ref{alx}),
(\ref{al}), (\ref{ko}) and algebras
$ {\cal B^{\alpha}}$ with
(and also without) the additional conditions (\ref{m2})
was developed by the present author in \cite{3},
where complete classification results were obtained
for the cases, when the number of components is not more than 4.

\bex \label{e1}
{\sl A Poisson bracket generated by the Lie algebra
of vector fields on an $n$-dimensional torus $T^n$.}

{\rm  In this case $N = n$. The commutator of
vector fields $\xi$ and $\eta$ has the form
\be
[\xi,\eta]_k=\xi_s\frac{\pa\eta_k}{\pa x^s}
-\eta_s\frac{\pa\xi_k}{\pa x^s}\,,
\ee
where $\xi = (\xi_1 (x), \ldots, \xi_n (x))$,
$\eta = (\eta_1 (x), \ldots, \eta_n (x))$,
$x = (x^1, \ldots, x^n)$, $x \in T^n$.

Then for the corresponding Hamiltonian operator
$M^{ij}$ we obtain
\bea
&&
\int_{T^n}u^k[\xi,\eta]_k\,d^nx
=\int_{T^n}u^k\left (\xi_i\frac{\pa\eta_k}{\pa x^i}
-\eta_j\frac{\pa\xi_k}{\pa x^j} \right )\,d^nx =\nn\\
&&
=\int_{T^n}\xi_i\biggr(u^i\frac d{d x^j}
+u^j\frac d{d x^i}
+\frac{\pa u^i}{\pa x^j}\biggr)\eta_j\,d^nx
=\int_{T^n}\xi_iM^{ij}\eta_j\,d^nx,
\label{c4}
\eea
\be
M^{ij}=u^i\frac d{d x^j}
+u^j\frac d{d x^i}+\frac{\pa u^i}{\pa x^j}
= \left ( u^i \delta^{j\alpha} +
u^j \delta^{i\alpha} \right )\frac d{d x^{\alpha}}
+ \delta^i_k \delta^{j\alpha}u^k_{\alpha}\,,
\label{c5}
\ee
\be
g^{ij\alpha}
=  u^i \delta^{j\alpha} +
u^j \delta^{i\alpha}, \ \ \ \
b^{ij\alpha}_k = \delta^i_k \delta^{j\alpha}.
\label{c5a}
\ee

It follows from the relations of Theorem
\ref{t1} that if for a $2$-cocycle on the Lie
algebra of vector fields on $T^n$ the corresponding Poisson
structure remains in the class of Poisson structures
of hydrodynamic type (\ref{1n})
 (see formulae (\ref{c0}) and (\ref{c3})), then
this 2-cocycle is cohomologous to zero.

For $n\le2$ the Poisson structure (\ref{c5}) is
nondegenerate, but for $n>2$ all the metrics in (\ref{c5})
are degenerate.

For $n = 1$ the only metric has the form $g (u) = 2 u$,
and the Hamiltonian operator has the form
\be
M=2 u\frac d{d x}
+ u_x.
\label{c5b}
\ee

For $n = 2$ both the metrics are nondegenerate
and indefinite:
\be
(g^{ij1})=\left ( \begin{array} {cc}
2 u^1&u^2\\ u^2&0\end{array} \right ),\qquad
(g^{ij2})=\left ( \begin{array} {cc}
 0&u^1\\ u^1&2u^2\end{array} \right ),
\ee
and the Hamiltonian operator has the form
\be
(M^{ij})=\left ( \begin{array} {cc}
2 u^1&u^2\\ u^2&0\end{array} \right ){d \over d x^1} +
\left ( \begin{array} {cc}
 0&u^1\\ u^1&2u^2\end{array} \right ) {d \over d x^2} +
\left ( \begin{array} {cc}
u^1_{x^1}&u^1_{x^2}\\ u^2_{x^1}&u^2_{x^2}\end{array} \right ).
\label{z1}
\ee

The Poisson structure (\ref{z1}) cannot be reduced
to a constant form by a local change of coordinate,
since the obstruction tensor $T^{ijk12}$ is not equal to zero
identically:
\bea
&&
T^{ijk\alpha\beta} =
g^{ks\beta} g^{ir\alpha} \left ( \Gamma^{j\beta}_{rs} -
\Gamma^{j\alpha}_{rs} \right ) = - g^{ir\alpha} b^{kj\beta}_r
+ g^{ks\beta} b^{ij\alpha}_s = \nn\\
&&
= - g^{ik\alpha}\delta^{j\beta}
+ g^{ki\beta} \delta^{j\alpha} =
- \left ( u^i \delta^{k\alpha}  +
u^k \delta^{i\alpha} \right ) \delta^{j\beta}
+ \left ( u^i \delta^{k\beta}
+ u^k \delta^{i\beta} \right ) \delta^{j\alpha} = \nn\\
&&
= u^i \left ( \delta^{j\alpha} \delta^{k\beta}  -
 \delta^{j\beta} \delta^{k\alpha} \right )
+  u^k \left ( \delta^{j\alpha} \delta^{i\beta} -
 \delta^{j\beta} \delta^{i\alpha} \right ).
\eea
In particular, $T^{11212} = u^1$.
}
\eex

This example is connected to the two-dimensional
Euler hydrodynamics of ideal incompressible
fluid (with further reduction to divergence-free
vector fields).

\bt {\rm \cite{3}}.
 \label{t4a}
If for $N=n=2$, for a nondegenerate Poisson structure
of the form {\rm (\ref{1n})}, the obstruction tensor
$T^{i12}_{jk} (u)$ is not equal to zero identically, i.e.,
this Poisson structure cannot be reduced to a constant form
by a local change of coordinates, then it
can be reduced to a canonical form generated
by the flat metrics
\be
(g^{ij1})=\left ( \begin{array} {cc}
1&0\\ 0&-1\end{array} \right ),\qquad
(g^{ij2})=\left ( \begin{array} {cc}
 2u^2&u^1+u^2\\ u^1+u^2&2u^1\end{array} \right ).
\label{c6}
\ee
Both the metrics {\rm (\ref{c6})} are
indefinite and therefore, in particular,
if one of the metrics of a two-dimensional
two-component nondegenerate
Poisson bracket of hydrodynamic type is positive or negative definite,
then this Poisson bracket can be reduced to a constant form.
The Poisson structure generated by the canonical flat metrics
{\rm (\ref{c6})} is connected to the Lie algebra
of vector fields on a two-dimensional torus $T^2$.
\et

By virtue of Theorem \ref{t3}, it is obvious that
the canonical Poisson structure generated by the metrics
(\ref{c6}) cannot be reduced to a constant form
by local changes of coordinates.
The obstruction tensor $T^{ijk12}$ of this Poisson structure
is not equal to zero:
$T^{ijk12} = g^{ks2} g^{ir1} (\Gamma^{j2}_{rs} -
\Gamma^{j1}_{rs}) = g^{ks2} \varepsilon^i
\delta^{ir} \Gamma^{j2}_{rs} = g^{ks2} \varepsilon^i
 \Gamma^{j2}_{is} = - \varepsilon^i
 b^{kj2}_i$, $\varepsilon^1 = 1$, $\varepsilon^2 = - 1$,
in particular,
$T^{21112}=b^{112}_2=1$.

The Poisson structure
$$
\{w^i(x),w^j(y)\}= \left (w^i(x)\delta^{j\alpha}+
w^j(x)\delta^{i\alpha} \right )\delta_{\alpha}(x-y)+
\delta^i_k \delta^{j\alpha}w^k_{\alpha}(x)\delta(x-y)
$$
generated by the Lie algebra of vector fields
on a two-dimensional torus $T^2$ (example \ref{e1})
is reduced to the canonical form
given by the metrics
{\rm (\ref{c6})} by the following
local quadratic change of coordinates \cite{3}:
\be
w^1=\frac12((u^1)^2-(u^2)^2),\ \ \
w^2=\frac12(u^1+u^2).\nn
\ee

Each two-dimensional two-component
nondegenerate Poisson structure of
hydrodynamic type either can be reduced to a constant
form or is generated by the Lie algebra
of vector fields on a two-dimensional torus
$T^2$ \cite{3}.

For $n>2$ each
multidimensional two-component
nondegenerate Poisson structure
of hydrodynamic type either can be reduced to a constant
form or can be reduced to the two-dimensional
canonical Poisson bracket
given by the metrics
{\rm (\ref{c6})}
by a local change of coordinates and by an unimodular change
of the independent space variables $x^i$ \cite{3}.

Note that each constant multidimensional (for an arbitrary dimension $n$)
two-component Dubrovin--Novikov bracket is reduced
by an unimodular change
of the independent space variables $x^i$ to a constant Dubrovin--Novikov
bracket whose dimension is not more than 3 (but it cannot be reduced,
generally speaking, to
a two-dimensional bracket).

\section{Compatible metrics}

Two pseudo-Riemannian contravariant
metrics $g^{ij}_1 (u)$ and
$g^{ij}_2 (u)$ are called {\it compatible} if
for any linear combination
$g^{ij} (u) = \lambda_1 g^{ij}_1 (u) + \lambda_2 g^{ij}_2 (u)$
of these metrics, where $\lambda_1$ and $\lambda_2$
are arbitrary constants, the coefficients of
the corresponding Levi-Civita connections
and the components of the corresponding Riemannian curvature
tensors are related by the same linear relation:
$\Gamma^{ij}_k (u)
= \lambda_1 \Gamma^{ij}_{1, k} (u) +
\lambda_2 \Gamma^{ij}_{2, k} (u),$
$R^{ij}_{kl} (u) = \lambda_1 R^{ij}_{1, kl} (u)
+ \lambda_2 R^{ij}_{2, kl} (u)$ \cite{5}.

Indices of the coefficients of connections and
indices of the curvature tensors are raised and lowered by
the metrics corresponding to them:
$\Gamma^{ij}_k (u)
= g^{is} (u) \Gamma^j_{sk} (u),$
$R^{ij}_{kl} (u) = g^{is} (u) R^j_{skl} (u)$.

If for any linear combination of metrics
the above-mentioned relation only for the corresponding Levi-Civita
connections is fulfilled, then the metrics are called
{\it almost compatible}.

The theory of compatible and almost compatible metrics
is constructed by the present author in the paper
\cite{5}. This theory is closely connected to
the theory of compatible Poisson brackets of
hydrodynamic type (local and nonlocal),
the theory of Frobenius manifolds and
integrable systems.
Recall that Poisson brackets are called
{\it compatible} if each their linear combination is a
Poisson bracket \cite{10}.

\bt {\rm \cite{6}}. \label{t8}
All metrics $g^{ij\alpha} (u)$, $1 \leq \alpha \leq n,$
defining a multidimensional Poisson bracket of the form
{\rm (\ref{1n})} are mutually compatible.
All one-dimensional Dubrovin--Novikov brackets forming
a multidimensional Poisson bracket
{\rm (}see Lemma {\rm \ref{l1})} are also mutually compatible.
\et

This important theorem immediately follows from
the tensor relations of
compatibility of metrics \cite{5}, \cite{11}.
For two arbitrary pseudo-Riemannian metrics
$g^{ij\alpha} (u)$ and $g^{ij\beta} (u)$ consider
an analog of the obstruction tensors that were defined above
for flat metrics of multidimensional Dubrovin--Novikov
brackets:
$T^{i\alpha\beta}_{jk} (u) =
\Gamma^{i\beta}_{jk} (u) - \Gamma^{i\alpha}_{jk} (u)$,
$T^{ijk\alpha\beta} (u) =
g^{ks\beta} (u) g^{ir\alpha} (u) T^{j\alpha\beta}_{rs} (u)$,
where
$\Gamma^{i\alpha}_{jk} (u)$ is the Levi-Civita connection
generated by the metric $g^{ij\alpha} (u)$.
Then we have the following theorem.

\bt  {\rm \cite{5}, \cite{11}}.
Any two pseudo-Riemannian
metrics $g^{ij\alpha} (u)$ and $g^{ij\beta} (u)$
are compatible if and only if the relations
{\rm (\ref{b1})} and {\rm (\ref{b3})}
of Theorem {\rm \ref{t2}} are fulfilled.
Moreover, the condition of almost compatibility
is equivalent to the relation {\rm (\ref{b1})}.
\et

Note that namely the relations (\ref{b1}) and (\ref{b3})
for flat metrics of multidimensional
Poisson brackets of hydrodynamic type  were found
by B.A.Dubrovin and S.P.Novikov in \cite{1}.

Thus the description and the classification of
multidimensional Poisson brackets of hydrodynamic type
correspond to the description and the classification of
an important special {\it subclass} of
compatible one-dimensional Poisson brackets of hydrodynamic type.
This subclass is singled out by the additional relations
(\ref{b2}) and (\ref{b4}).

The problem of the description of
compatible one-dimensional Poisson brackets of
hydrodynamic type that is equivalent to
the description of {\it flat pencils of
compatible metrics} or, in other words,
{\it local quasi-Frobenius manifolds}, plays an important role
in the theory of integrable systems of hydrodynamic type,
the theory of Frobenius manifolds, the modern differential geometry
and mathematical physics
(see \cite{12}).

The problem of the description of all nonsingular
pairs of compatible flat metrics was solved by the present author
in \cite{13}, \cite{14} (see also \cite{11}), where
were obtained the nonlinear equations
describing all nonsingular pairs of
compatible flat metrics (these nonlinear equations are
a nontrivial nonlinear reduction of
the Lam\'e equations defining all the curvilinear
orthogonal coordinate systems in a pseudo-Euclidean space)
and was proved the integrability
of these nonlinear equations by the inverse scattering method, moreover,
was found an explicit integration and linearization procedure, which
reduces the integration of these nonlinear equations
to the solution of linear problems
(thus a local classification of semisimple
quasi-Frobenius manifolds was also obtained); see also \cite{5},
where the results were announced. Afterwards for these nonlinear equations
also a Lax pair was found in \cite{14a}.

A pair of metrics $g^{ij}_1 (u)$ and $g^{ij}_2 (u)$
is called {\it nonsingular} if all eigenvalues of
this pair of metrics, i.e., the roots of the equation
$\det (g^{ij}_1 (u) - \lambda g^{ij}_2 (u)) = 0,$
are distinct (the situation of a pair of metrics in general position).

\bt {\rm \cite{5}}. \label{t9}
If a pair of metrics $g^{ij}_1 (u)$ and $g^{ij}_2 (u)$
is nonsingular, then these metrics are compatible
if and only if for the affinor
$v^i_j (u) = g^{is}_1 (u) g_{2, sj} (u)$
the Nijenhuis tensor
$$N^k_{ij} (u) = v^s_i (u) {\pa v^k_j \over \pa u^s}
-  v^s_j (u) {\pa v^k_i \over \pa u^s} +
  v^k_s (u) {\pa v^s_i \over \pa u^j} -
  v^k_s (u) {\pa v^s_j \over \pa u^i}$$
vanishes
{\rm (}here $g_{2, sj} (u)$ is the corresponding covariant metric,
$g^{is}_2 (u) g_{2, sj} (u) = \delta^i_j${\rm )}.
\et

In this case, by virtue of the Nijenhuis theorem \cite{15}
there exist local coordinates (generally speaking, complex)
in which the affinor $v^i_j (u)$ is diagonal in a domain:
$v^i_j (u) = \lambda^i (u) \delta^i_j.$
Then, in these special local coordinates, both the metrics
$g^{ij}_1 (u)$ and $g^{ij}_2 (u)$ are also necessarily
diagonal. Indeed, the eigenvalues $\lambda^i (u)$,
$1 \leq i \leq N,$ of the affinor $v^i_j (u)$ coincide with
the eigenvalues of the pair of metrics $g^{ij}_1 (u)$
and $g^{ij}_2 (u)$ and, by our condition, they are distinct
(the pair of metrics is nonsingular):
$\lambda^i (u) \neq \lambda^j (u)$ for $i \neq j$.
Then, in the special local coordinates, in which the affinor
$v^i_j (u)$ is
diagonal in a domain, we have
$g^{ij}_1 (u) = \lambda^i (u) g^{ij}_2 (u)$.
By virtue of symmetry of metrics we obtain
$(\lambda^i (u) - \lambda^j (u)) g^{ij}_2 (u) = 0$,
therefore if $i \neq j$, then $g^{ij}_1 (u) = g^{ij}_2 (u) = 0$,
i.e., the metrics are diagonal.
Moreover, the following important theorem is valid:

\bt {\rm \cite{5}}. \label{t10}
If a pair of metrics $g^{ij}_1 (u)$ and $g^{ij}_2 (u)$ is
nonsingular, then these metrics are compatible
if and only if there exist local coordinates
$u = (u^1, \ldots, u^N)$
(possibly, complex) such that
$g^{ij}_2 (u) = g^i (u) \delta^{ij}$
and $g^{ij}_1 (u) = f^i (u^i) g^i (u) \delta^{ij}$,
where $f^i (u^i),$
$i = 1, \ldots, N,$ are functions of one variable
(possibly, complex).
\et

Note here that for an arbitrary pair of metrics
the condition of vanishing the corresponding Nijenhuis tensor
is equivalent to the condition of almost compatibility of the metrics,
but not compatibility (there are
corresponding counterexamples) \cite{5}.

\section{Classification theorem}

\bt  \label{t11}
If for a nondegenerate multidimensional Poisson bracket
of the form {\rm (\ref{1n})} one of the metrics
$g^{ij\alpha} (u)$ forms nonsingular pairs with all
the remaining metrics of the bracket, then
this Poisson bracket can be reduced to a constant form
by a local change of coordinates.
\et

Let us prove that all the obstruction tensors
$T^{i\alpha\beta}_{jk} (u)$ are equal to zero identically.
Without loss of generality we can consider that
the metric
$g^{ij1} (u)$ forms nonsingular pairs with all the remaining metrics
of the bracket. Let $\beta \neq 1$.
By theorem \ref{t8} the metrics $g^{ij\beta} (u)$ and
$g^{ij1} (u)$ are compatible.
By Theorem \ref{t10} there exist local coordinates
$u = (u^1, \ldots, u^N)$
(possibly, complex) such that $g^{ij\beta} (u) = g^{i} (u) \delta^{ij}$
and $g^{ij1} (u) = f^{i} (u^i)
g^{i} (u) \delta^{ij}$, where $f^i (u^i),$
$i = 1, \ldots, N,$ are functions of one variable
(possibly, complex). In these local coordinates
$\Gamma^{i\beta}_{jk} (u) = 0$ if all the indices
$i, j, k$ are distinct;
$\Gamma^{i\beta}_{ik} (u) = \Gamma^{i\beta}_{ki} (u) =
- (1/(2g^i (u)))(\pa g^i/\pa u^k)$ for any indices $i, k$;
$\Gamma^{i\beta}_{jj} (u) =
(g^i (u)/(2(g^j (u))^2))(\pa g^j/\pa u^i)$ for $i \neq j$.
For the obstruction tensor
$T^{ijk1\beta} (u)$ in the local coordinates under consideration
we obtain: $T^{ijk1\beta} (u) = 0$ if all the indices
$i, j, k$ are distinct; $T^{ijj1\beta} (u) =
T^{iij1\beta} (u) = 0$ for $i \neq j$;
$T^{iii1\beta} (u) = ((g^i (u))^2/2)(d f^i (u^i) /d u^i)$;
$T^{iji1\beta} (u) =
((f^i (u^i) - f^j (u^j))/2)(\pa g^i/\pa u^j)$ for $i \neq j$.
From the relation (\ref{b2}) of Theorem \ref{t2}
we obtain that $T^{iii1\beta} = 0$,
$f^i (u^i) = \mu^i = {\rm const}$, and $T^{iji1\beta} (u) = 0$
for $i \neq j$, $g^i (u) = g^i (u^i)$ are functions of one variable.
Thus $T^{ijk1\beta} (u) = 0$ for any $\beta$.
Since for all indices $\alpha$ and $\beta$
the identity $T^{i\alpha\beta}_{jk} (u) =
T^{i1\beta}_{jk} (u) - T^{i1\alpha}_{jk} (u)$ holds,
we obtain that $T^{i\alpha\beta}_{jk} (u) = 0$ for all
$\alpha$ and $\beta$.

Thus all pairs of metrics defining Dubrovin--Novikov
brackets that cannot be reduced to a constant form,
in particular, all nonconstant canonical pairs of metrics
(they are connected to nontrivial Lie algebras of
hydrodynamic type and nontrivial quasi-Frobenius
algebras), are {\it nonsingular}, i.e., they have
coinciding eigenvalues. In particular,
for the unique nonconstant canonical pair
of two-component metrics
(\ref{c6}) we have: $\det (g^{ij2} - \lambda g^{ij1}) = 0,$
$\lambda_1 = \lambda_2 = u^2 - u^1$.
The theory of singular nondegenerate multidimensional (and also all
degenerate) Dubrovin--Novikov brackets and algebras related to them
is far from completeness at present.

{\bf Acknowledgements.} I gratefully acknowledge the hospitality of
Max-Planck-Institut f\"ur Mathematik in Bonn (Germany).

%\begin{flushleft}
%e-mail: mokhov@mi.ras.ru; mokhov@landau.ac.ru; mokhov@bk.ru\\
%\end{flushleft}

\end{document}